\documentclass[a4paper,12pt]{article}
\usepackage[russian]{babel}
\usepackage{ucs}
\usepackage{amsfonts, amssymb, amsmath, amsthm, amscd}
\usepackage{cite}
\ifx\undefined \pdfgentounicode \else
\input{glyphtounicode} \pdfgentounicode=1
\fi

\author{A.A. Vasil'eva}
\title{Kolmogorov widths of a finite‐dimensional balls intersection}
\date{}
\begin{document}

\maketitle

\newenvironment{Biblio}{%
                  \renewcommand{\refname}{\footnotesize REFERENCES}%
                  \begin{thebibliography}{99}%
                  \itemsep -0.2em}%
                  {\end{thebibliography}}

\def\inff{\mathop{\smash\inf\vphantom\sup}}
\renewcommand{\le}{\leqslant}
\renewcommand{\ge}{\geqslant}
\newcommand{\sgn}{\mathrm {sgn}\,}
\newcommand{\inter}{\mathrm {int}\,}
\newcommand{\dist}{\mathrm {dist}}
\newcommand{\supp}{\mathrm {supp}\,}
\newcommand{\R}{\mathbb{R}}
\renewcommand{\C}{\mathbb{C}}
\newcommand{\Z}{\mathbb{Z}}
\newcommand{\N}{\mathbb{N}}
\newcommand{\Q}{\mathbb{Q}}
\theoremstyle{plain}
\newtheorem{Trm}{Theorem}
\newtheorem{trma}{Theorem}

\newtheorem{Def}{Definition}
\newtheorem{Cor}{Corollary}
\newtheorem{Lem}{Lemma}
\newtheorem{Rem}{Remark}
\newtheorem{Sta}{Proposition}
\newtheorem{Sup}{Assumption}
\newtheorem{Supp}{Assumption}
\newtheorem{Exa}{Example}
\renewcommand{\proofname}{\bf Proof}
\renewcommand{\thetrma}{\Alph{trma}}
\renewcommand{\theSupp}{\Alph{Supp}}

\section{Introduction}

When studying the problem of estimating the widths of Sobolev weight classes with restrictions on derivatives in various $L_p $-metrics \cite{vas_1},  the question arose about the widths of intersections of two finite-dimensional balls. Their order estimates were obtained in \cite{vas_finite}. Here we generalize this result to the intersection of an arbitrary family of balls. 

Let $X$ be a normed space, $M\subset X$, $n\in \Z_+$. We denote by ${\cal L}_n(X)$ the family of subspaces in $X$
of dimension at most $n$. The Kolmogorov $n$-width of the set $M$ in the space $X$ is defined as follows:
$$
d_n(M, \, X) = \inf _{L\in {\cal L}_n(X)} \sup _{x\in M} \inf
_{y\in L} \|x-y\|.
$$

Given $N\in \N$, $1\le p\le \infty$, we denote by $l_p^N$ the space $\R^N$ with the norm
$$
\|(x_1, \, \dots, \, x_N)\|_p = \left\{ \begin{array}{l} \left(\sum \limits _{i=1}^N |x_i|^p\right)^{1/p}, \quad \text{if } 1\le p<\infty, \\ \max _{1\le i\le N} |x_i|, \quad \text{if }p=\infty.\end{array}\right.
$$
By $B_p^N$ we denote the unit ball in the space $l_p^N$.

Let $A$ be a non‐empty set, $1\le p_\alpha \le \infty$, $\nu_\alpha>0$ for every $\alpha \in A$ and $p_\alpha \ne p_\beta$ for $\alpha \ne \beta$, and let
\begin{align}
\label{m_set_def}
M = \cap _{\alpha \in A} \nu _\alpha B_{p_\alpha}^N. 
\end{align}

In this paper we obtain order estimates of the widths $d_n(M, \, l_q^N)$ for $n\le N/2$.

The explicit values of $d_n(B_p^N, \, l_q^N)$ were obtained by Pietsch and Stesin \cite{pietsch1, stesin} for $p\ge q$, and by Kolmogorov, Petrov, Smirnov \cite{k_p_s} and Stechkin \cite{stech_poper} for $p=1$, $q=2$. For $p<q<\infty$, Gluskin \cite{gluskin1, bib_gluskin} obtained order estimates of $d_n(B_p^N, \, l_q^N)$. The case $q=\infty$ was investigated by Kashin \cite{kashin_oct, bib_kashin}, Gluskin and Garnaev \cite{garn_glus}; for $p\ge 2$, order estimates were obtained, and for $1\le p<2$, the upper and the lower estimates differ by a logarithmic multiplier. For details, see \cite{itogi_nt}.

Galeev \cite{galeev1} obtained order estimates for
$d_n(M, \, l_q^{2n})$, where $M$ is defined by (\ref{m_set_def}), $A\subset \R$, $\nu_\alpha = (2n)^{‐\alpha}$. This result was applied for estimating the widths of intersections of some function classes \cite{galeev1, galeev4}.

Given $\alpha$, $\beta\in A$, we set
\begin{align}
\label{kappa}
\varkappa _{\alpha, \, \beta} = \left\{ \begin{array}{l}\left(\frac{\nu_\beta}{\nu_\alpha}\right)^{\frac{p_\alpha p_\beta}{p_\alpha ‐p_\beta}}, \quad \text{if }\alpha\ne \beta, \\ 1, \quad \text{if }\alpha=\beta.\end{array}\right.
\end{align}
Then
\begin{align}
\label{kappa_ab} \nu _\alpha \varkappa _{\alpha,\beta}^{‐1/p_\alpha} = \nu_\beta \varkappa _{\alpha,\beta} ^{‐1/p_\beta}.
\end{align}
Notice that $\varkappa _{\alpha,\beta} = \varkappa _{\beta,\alpha}$.

Suppose that 
\begin{align}
\label{kappa_n}
1\le \varkappa _{\alpha,\beta} \le N, \quad \alpha, \, \beta\in A
\end{align}
(in \S 3 it will be shown how to reduce the general case to this case).
Then $\{\nu_\alpha:\alpha \in A\}$ is bounded and separated from zero.

Now we formulate the main result of this paper. First we give notation for order equalities and inequalities.

Let $X$, $Y$ be sets, and let $f_1$, $f_2:\ X\times Y\rightarrow \mathbb{R}_+$.
We write $f_1(x, \, y)\underset{y}{\lesssim} f_2(x, \, y)$ (or
$f_2(x, \, y)\underset{y}{\gtrsim} f_1(x, \, y)$) if for any
$y\in Y$ there exists $c(y)>0$ such that $f_1(x, \, y)\le
c(y)f_2(x, \, y)$ for every $x\in X$; $f_1(x, \,
y)\underset{y}{\asymp} f_2(x, \, y)$ if $f_1(x, \, y)
\underset{y}{\lesssim} f_2(x, \, y)$ and $f_2(x, \,
y)\underset{y}{\lesssim} f_1(x, \, y)$.

We denote 
\begin{align}
\label{k_def}
K=\{(\nu_\alpha, \, 1/p_\alpha)\}_{\alpha \in A}.
\end{align}
\begin{Trm}
\label{main}
Let $n, \, N\in \N$, $n\le N/2$, and let $M$ be defined by (\ref{m_set_def}). Suppose that (\ref{kappa_n}) holds.
\begin{enumerate}
\item Let $p_\alpha \ge q$ for every $\alpha \in A$. Then
\begin{align}
\label{1}
d_n(M, \, l_q^N) \asymp \inf _{\alpha \in A} \nu_\alpha N^{1/q‐1/p_\alpha}.
\end{align}

\item Let $q\le 2$, $p_\alpha \le q$ for any $\alpha \in A$. Then 
\begin{align}
\label{2}
d_n(M, \, l_q^N) \asymp \inf _{\alpha \in A} \nu_\alpha.
\end{align}

\item Let $q\le 2$, $A_1=\{\alpha \in A:\; p_\alpha>q\}\ne \varnothing$, $A_2=\{\alpha \in A:\; p_\alpha<q\}\ne \varnothing$, $A_1' = \{\alpha \in A:\; p_\alpha\ge q\}$, $A_2' = \{\alpha \in A:\; p_\alpha\le q\}$. Then
\begin{align}
\label{3}
d_n(M, \, l_q^N) \asymp \inf \{\nu_\alpha \varkappa _{\alpha,\beta}^{1/q‐1/p_\alpha}:\; \alpha \in A_1' , \; \beta \in A_2' \}.
\end{align}

\item Let $2<q<\infty$, and let $p_\alpha\le 2$ for every $\alpha \in A$. Then
\begin{align}
\label{4}
d_n(M, \, l_q^N) \underset{q}{\asymp} \inf _{\alpha \in A} \nu_\alpha \min \{1, \, n^{‐1/2}N^{1/q}\}.
\end{align}
\item Let $2<q<\infty$, $A_1' =\{\alpha \in A:\; p_\alpha\ge q\}$, $A_2' = \{\alpha \in A:\; 2\le p_\alpha\le q\}$, $A_3' = \{\alpha \in A:\; p_\alpha \le 2\}$. Suppose that $A \not \subset A_1'$, $A \not \subset A_3'$. Let $K$ be defined by (\ref{k_def}). Then 
\begin{align}
\label{5}
\begin{array}{c} d_n(M, \, l_q^N) \underset{q}{\asymp} \Phi(n, \, N, \, q, \, K):= \\ =\min \{\Phi_1(n, \, N, \, q, \, K), \, \Phi_2(n, \, N, \, q, \, K), \, \Phi_3(n, \, N, \, q, \, K)\}, \end{array}
\end{align}
where
\begin{align}
\label{6}
\Phi_1(n, \, N, \, q, \, K) =\inf \{\nu_\alpha \varkappa _{\alpha,\beta}^{1/q‐1/p_\alpha}:\; \alpha \in A_1', \; \beta \in A_2'\cup A_3'\},
\end{align}
\begin{align}
\label{7}
\Phi_2(n, \, N, \, q, \, K) =\inf \left \{ \nu_\alpha (\min \{1, \, n^{‐1/2}N^{1/q}\})^{\frac{1/p_\alpha‐1/q}{1/2‐1/q}}, \; \alpha \in A_2'\right\},
\end{align}
\begin{align}
\label{8}
\Phi_3(n, \, N, \, q, \, K) = \min \{\nu_\alpha \varkappa _{\alpha,\beta}^{1/2‐1/p_\alpha}\min \{1, \, n^{‐1/2}N^{1/q}\}:\; \alpha \in A_1'\cup A_2', \; \beta \in A_3'\}
\end{align}
(the infimum of the empty set is $+\infty$).
\end{enumerate}
\end{Trm}

\begin{Rem}
If $p_\alpha=q$, we have
$$
\nu_\alpha \varkappa _{\alpha,\beta}^{1/q‐1/p_\alpha} = \nu_\alpha = \nu_\alpha (\min \{1, \, n^{‐1/2}N^{1/q}\})^{\frac{1/p_\alpha‐1/q}{1/2‐1/q}};
$$
if $p_\beta = q$, then
$$
\nu_\alpha \varkappa _{\alpha,\beta}^{1/q‐1/p_\alpha} \stackrel{(\ref{kappa_ab})}{=}\nu_\beta \varkappa _{\alpha,\beta}^{1/q‐1/p_\beta} =\nu_\beta=\nu_\beta (\min \{1, \, n^{‐1/2}N^{1/q}\})^{\frac{1/p_\beta‐1/q}{1/2‐1/q}}.
$$
If $p_\alpha=2$, we have
$$
\nu_\alpha (\min \{1, \, n^{‐1/2}N^{1/q}\})^{\frac{1/p_\alpha‐1/q}{1/2‐1/q}} = $$$$= \nu_\alpha \min \{1, \, n^{‐1/2}N^{1/q}\} =\nu_\alpha \varkappa _{\alpha,\beta}^{1/2‐1/p_\alpha}\min \{1, \, n^{‐1/2}N^{1/q}\};
$$
if $p_\beta=2$, then
$$
\nu_\alpha \varkappa _{\alpha,\beta}^{1/2‐1/p_\alpha}\min \{1, \, n^{‐1/2}N^{1/q}\} \stackrel{(\ref{kappa_ab})}{=} \nu_\beta \varkappa _{\alpha,\beta}^{1/2‐1/p_\beta}\min \{1, \, n^{‐1/2}N^{1/q}\}=$$$$ = \nu_\beta \min \{1, \, n^{‐1/2}N^{1/q}\}=\nu_\beta (\min \{1, \, n^{‐1/2}N^{1/q}\})^{\frac{1/p_\beta‐1/q}{1/2‐1/q}}.
$$
\end{Rem}

\begin{Trm}
\label{main_infty} Let $n, \, N\in \N$, $n\le N‐1$, let $M$ be defined by (\ref{m_set_def}), $p_\alpha\ge 2$ for every $\alpha \in A$. Suppose that (\ref{kappa_n}) holds. Then
\begin{align}
\label{d_n_infty_trm}
d_n(M, \, l_\infty^N) \asymp \inf \left\{ \nu _\alpha \min \left(1, \, n^{‐1} \log \frac{2N}{n}\right)^{1/p_\alpha}:\; \alpha \in A\right\}.
\end{align}
\end{Trm}

\section{Proof of the main results}

First we formulate some well‐known results.
\begin{trma}
\label{glus} {\rm \cite{bib_gluskin}} Let $1\le p\le q<\infty$,
$0\le n\le N/2$.
\begin{enumerate}
\item Let $1\le q\le 2$. Then $d_n(B_p^N, \, l_q^N) \underset{p,q}{\asymp}
1$.

\item Let $2<q<\infty$, $\lambda_{pq} =\min \left\{1, \,
\frac{1/p-1/q}{1/2-1/q}\right\}$. Then
$$
d_n(B_p^N, \, l_q^N) \underset{p,q}{\asymp} \min \{1, \,
n^{-1/2}N^{1/q}\} ^{\lambda_{pq}}.
$$
\end{enumerate}
\end{trma}

\begin{trma}
\label{p_s} {\rm \cite{pietsch1, stesin}} Let $1\le q\le p\le
\infty$, $0\le n\le N$. Then
\begin{align}
\label{ps1}
d_n(B_p^N, \, l_q^N) = (N-n)^{1/q-1/p}.
\end{align}
\end{trma}

\begin{trma}
\label{garn_glusk_teor} {\rm \cite{garn_glus}}
Let $p\ge 2$, $n\le N‐1$. Then
\begin{align}
\label{glusk_dn2} d_n(B_p^N, \, l_\infty^N) \asymp
\min\{1, \, n^{-1/p}\log^{1/p}(1+N/n)\}.
\end{align}
\end{trma}

Given $k\in \{1,\, \dots, \, N\}$, we define the set $V_k\subset \R^N$ by
$$
V_k={\rm conv}\, \{(\varepsilon_1 \hat{x}_{\sigma(1)}, \, \dots,
\, \varepsilon_N \hat{x}_{\sigma(N)}):\; \varepsilon_j=\pm 1,
\;1\le j\le N, \; \sigma \in S_N\},
$$
where $\hat{x}_j=1$ for $1\le j\le k$, $\hat{x}_j=0$ for $k+1\le j\le
N$, $S_N$ is the group of permutations of $N$ elements. Notice that $V_1
= B_1^N$, $V_N = B_\infty^N$.

For $2\le q<\infty$, the lower estimates for $d_n(V_k, \, l_q^N)$ were obtained by Gluskin \cite{gluskin1}.

\begin{trma}
\label{gl_q_g_2} {\rm \cite{gluskin1}} Let $2\le q<\infty$, $1\le k\le N$, $n\le 
\min \{N^{\frac 2q}k^{1 -\frac 2q}, \, N/2\}$. Then
\begin{align}
\label{kq1} d_n(V_k, \, l_q^N) \underset{q}{\gtrsim} k^{1/q}.
\end{align}
Let $2\le q<\infty$, $1\le k\le N$, $N^{\frac 2q}k^{1 -\frac 2q}\le n\le N/2$. Then
\begin{align}
\label{kq0} d_n(V_k, \, l_q^N) \underset{q}{\gtrsim}
k^{1/2}n^{-1/2}N^{1/q}.
\end{align}
\end{trma}

The following result was obtained by Gluskin \cite{bib_glus_3} (the constant in the order inequality depends on $q$), and by Malykhin and Rjutin \cite{mal_rjut} (the constant does not depend on $q$). In \cite[p. 39]{bib_glus_3} it was noticed that Galeev obtained the equality $d_n(V_k, \, l_1^N) =\min\{k, \, N‐n\}$.
\begin{trma}
 {\rm \cite{bib_glus_3, mal_rjut}.} Let $1\le q\le 2$, $n\le
N/2$. Then 
\begin{align}
\label{gl_q_l_2}
d_n(V_k, \, l_q^N) \gtrsim k^{1/q}.
\end{align}
\end{trma}

The following result was obtained by Konyagin, Malykhin and Rjutin \cite{kon_m_r}.
\begin{trma} \label{t_kmr}
{\rm \cite{kon_m_r}.} Let $n\le N‐1$, $1\le k\le N$. Then there is an absolute constant $C\ge 1$ such that for $k\ge C n\log ^{‐1}\frac {2N}{n}$ $$d_n(V_k, \, l_\infty^N) \ge 1/2.$$
\end{trma}
\renewcommand{\proofname}{\bf Proof of Theorem \ref{main}}
\begin{proof}
Let $K$ be defined by (\ref{k_def}). From (\ref{kappa_n}) it follows that $K$ is bounded and $\inf _{\alpha \in A}\nu_\alpha >0$. We show that it suffices to consider the case of a compact $K$. Indeed, let $\overline{K}=\{(\nu_\alpha, \, 1/p_\alpha):\; \alpha \in A'\}$ be the closure of $K$. For $\alpha$, $\beta\in A'$ we define the numbers $\varkappa _{\alpha,\beta}$ by (\ref{kappa}). We prove that
\begin{enumerate}
\item $M = \cap _{\alpha\in A'} \nu_\alpha B_{p_\alpha}^N$,
\item on the right‐hand sides of (\ref{1})‐‐(\ref{8}) the set $A$ can be replaced by $A'$, and the minimum is attained.
\end{enumerate}

Let $\nu_{\alpha_k} \underset{k\to \infty}{\to} \nu_\alpha$, $1/p_{\alpha_k}\underset{k\to \infty}{\to} 1/p_\alpha$. Then $\cap _{k\in \N} \nu_{\alpha_k} B_{p_{\alpha_k}}^N \subset \nu_\alpha B_{p_\alpha}^N$. This implies the first assertion.

Let $(\nu_{\alpha_k}, 1/p_{\alpha_k}) \underset{k\to \infty}{\to}(\nu_\alpha, \, 1/p_\alpha)$, $(\nu_{\beta_k}, 1/p_{\beta_k}) \underset{k\to \infty}{\to}(\nu_\beta, \, 1/p_\beta)$, $\alpha\ne \beta$. Then
\begin{align}
\label{lim_kappa} \varkappa _{\alpha_k, \, \beta_k} \underset{k\to\infty}{\to} \varkappa _{\alpha,\beta}.
\end{align}
In particular, $1\le \varkappa _{\alpha,\beta}\le N$ for all $\alpha$, $\beta\in A'$.

It is clear that if $(\nu_{\alpha_k}, 1/p_{\alpha_k}) \underset{k\to \infty}{\to}(\nu_\alpha, \, 1/p_\alpha)$, then $\nu_{\alpha_k} \underset{k\to \infty}{\to} \nu_\alpha$, $\nu _{\alpha_k}N^{1/q‐1/p_{\alpha_k}}\underset{k\to \infty}{\to} \nu_\alpha N^{1/q‐1/p_\alpha}$; if, in addition, $q>2$, then $$\nu _{\alpha_k}(\min \{1, \, n^{‐1/2}N^{1/q}\}) ^{\frac{1/p_{\alpha_k}‐1/q}{1/2‐1/q}} \underset{k\to \infty}{\to} \nu_\alpha (\min \{1, \, n^{‐1/2}N^{1/q}\}) ^{\frac{1/p_{\alpha}‐1/q}{1/2‐1/q}}.$$

Let $(\nu_{\alpha_k}, 1/p_{\alpha_k}) \underset{k\to \infty}{\to}(\nu_\alpha, \, 1/p_\alpha)$, $(\nu_{\beta_k}, 1/p_{\beta_k}) \underset{k\to \infty}{\to}(\nu_\beta, \, 1/p_\beta)$, $p_{\alpha_k}\ge s$, $p_{\beta_k}\le s$ for every $k$. We show that $\nu_{\alpha_k} \varkappa _{\alpha_k,\beta_k}^{1/s‐1/p_{\alpha_k}} \underset{k\to \infty}{\to} \nu_\alpha \varkappa _{\alpha,\beta} ^{1/s‐1/p_\alpha}$. If $\alpha\ne \beta$, it follows from (\ref{lim_kappa}). If $\alpha=\beta$, then $p_\alpha=p_\beta=s$, $\varkappa _{\alpha,\beta}=1$, $1\le \varkappa _{\alpha_k,\beta_k}^{1/s‐1/p_{\alpha_k}}\le N^{1/s‐1/p_{\alpha_k}}\underset{k\to \infty}{\to} 1$.

This completes the proof of the second assertion.

In what follows we suppose that $K$ is compact. 

Let $p_\beta\le s\le p_\alpha$. Then $\frac{1}{s}= \frac{1‐\lambda}{p_\beta} +\frac{\lambda}{p_\alpha}$ for some $\lambda \in [0, \, 1]$. From H\"{o}lder's inequality or from \cite[Theorem 2]{galeev1} it follows that
\begin{align}
\label{s_p_alpha} B_{p_\alpha}^N \cap \varkappa _{\alpha,\beta}^{1/p_\beta‐1/p_\alpha} B_{p_\beta}^N \subset \varkappa _{\alpha,\beta}^{(1‐\lambda)(1/p_\beta‐1/p_\alpha)} B_s^N = \varkappa _{\alpha,\beta}^{1/s‐1/p_\alpha} B_s^N.
\end{align}

Now we prove (\ref{1})‐‐(\ref{5}).

{\bf Case 1.} The upper estimate holds, since for every $\alpha \in A$
$$
d_n(M, \, l_q^N) \le \nu_\alpha d_n(B_{p_\alpha}^N, \, l_q^N)\stackrel{(\ref{ps1})}{=}\nu_\alpha(N‐n)^{1/q‐1/p_\alpha} \le \nu_\alpha N ^{1/q‐1/p_\alpha}.
$$
Let us prove the lower estimate. Let 
\begin{align}
\label{nu_a_st_n_qpa}
\nu_{\alpha_*} N ^{1/q‐1/p_{\alpha_*}} = \min _{\alpha \in A} \nu_\alpha N ^{1/q‐1/p_\alpha}. 
\end{align}
We prove that 
\begin{align}
\label{nu_a_n_p_bn_m}
\nu_{\alpha_*} N^{‐1/p_{\alpha_*}}B^N_\infty \subset M. 
\end{align}
It suffices to check that
$\nu_{\alpha_*}N^{1/p_\alpha ‐1/p_{\alpha_*}}\le \nu _{\alpha}$ for all $\alpha\in A$. Indeed, this inequality holds, since
$\nu_{\alpha_*} N ^{1/q‐1/p_{\alpha_*}} \le \nu_\alpha N ^{1/q‐1/p_\alpha}$ by (\ref{nu_a_st_n_qpa}).

Hence,
$$
d_n(M, \, l_q^N) \stackrel{(\ref{nu_a_n_p_bn_m})}{\ge} \nu_{\alpha_*}N^{‐1/p_{\alpha_*}} d_n(B^N_\infty, \, l_q^N) \stackrel{(\ref{ps1})}{=} \nu_{\alpha_*}N^{‐1/p_{\alpha_*}} (N‐n)^{1/q} \gtrsim \nu_{\alpha_*}N^{1/q‐1/p_{\alpha_*}}
$$
(recall that $n\le N/2$).

{\bf Cases 2 and 4.} Let $\nu _{\alpha_*}= \min _{\alpha \in A}\nu_\alpha$. Then $\nu _{\alpha_*}B_1^N \subset M \subset \nu _{\alpha_*}B^N_{\min\{q, \, 2\}}$. It remains to apply Theorem \ref{glus}.

{\bf Case 3.} Let us prove the upper estimate. Let $\alpha \in A_1'$, $\beta \in A_2'$. Then $$\nu _\alpha B_{p_\alpha}^N\cap \nu_\beta B_{p_\beta}^N \stackrel{(\ref{kappa_ab})}{=} \nu _\alpha (B_{p_\alpha}^N\cap \varkappa _{\alpha,\beta}^{1/p_\beta‐1/p_\alpha} B_{p_\beta}^N) \stackrel{(\ref{s_p_alpha})}{\subset} \nu_\alpha \varkappa _{\alpha,\beta}^{1/q‐1/p_\alpha}B_q^N.$$ Hence,
$$
d_n(M, \, l_q^N) \le \nu_\alpha \varkappa _{\alpha,\beta}^{1/q‐1/p_\alpha}d_n(B_q^N, \, l_q^N) = \nu_\alpha \varkappa _{\alpha,\beta}^{1/q‐1/p_\alpha}.
$$
It remains to take the minimum over $\alpha\in A_1'$, $\beta \in A_2'$.

Now we obtain the lower estimate. Let 
\begin{align}
\label{nu_a_min}
\nu_{\alpha_*} \varkappa _{\alpha_*,\beta_*}^{1/q‐1/p_{\alpha_*}} = \min _{\alpha \in A_1',\beta \in A_2'} \nu_\alpha \varkappa _{\alpha,\beta}^{1/q‐1/p_\alpha}.
\end{align}

First we consider the case $p_{\alpha_*}>q$, $p_{\beta_*}<q$.

Let $k = [\varkappa _{\alpha_*,\beta_*}]$. By (\ref{kappa_n}), we have $1\le k\le N$. We prove that 
\begin{align}
\label{nu_k_m}
\nu _{\alpha_*}k^{‐1/p_{\alpha_*}}V_k \subset 2M.  
\end{align}
It suffices to check that
\begin{align}
\label{nu_a_g}
\nu _{\alpha_*}\varkappa _{\alpha_*,\beta_*}^{1/p_\gamma‐1/p_{\alpha_*}}\le \nu _\gamma, \quad \gamma \in A
\end{align}
(then $\nu _{\alpha_*}k^{1/p_\gamma‐1/p_{\alpha_*}}\le 2\nu _\gamma$ for all $\gamma \in A$, and (\ref{nu_k_m}) holds).

Let $\gamma \in A_2'$. We have
$$
\nu_{\alpha_*} \varkappa _{\alpha_*,\beta_*}^{1/q‐1/p_{\alpha_*}} \stackrel{(\ref{nu_a_min})}{\le} \nu_{\alpha_*} \varkappa _{\alpha_*,\gamma}^{1/q‐1/p_{\alpha_*}};
$$
since $p_{\alpha_*}>q$, the inequality $\varkappa _{\alpha_*,\beta_*} \le \varkappa _{\alpha_*,\gamma}$ holds. By the conditions $p_{\alpha_*} >q\ge p_\gamma$, we get $\nu _{\alpha_*}\varkappa _{\alpha_*,\beta_*}^{1/p_\gamma‐1/p_{\alpha_*}}\le \nu _{\alpha_*}\varkappa _{\alpha_*,\gamma}^{1/p_\gamma‐1/p_{\alpha_*}} \stackrel{(\ref{kappa_ab})}{=} \nu_\gamma$.

Let $\gamma \in A_1'$. We have
$$
\nu_{\alpha_*} \varkappa _{\alpha_*,\beta_*}^{1/q‐1/p_{\alpha_*}} \stackrel{(\ref{nu_a_min})}{\le} \nu_{\gamma} \varkappa _{\gamma,\beta_*}^{1/q‐1/p_{\gamma}}.
$$
By (\ref{kappa_ab}), it is equivalent to 
$$
\nu_{\beta_*} \varkappa _{\alpha_*,\beta_*}^{1/q‐1/p_{\beta_*}} \le \nu_{\beta_*} \varkappa _{\gamma,\beta_*}^{1/q‐1/p_{\beta_*}}.
$$
Since $p_{\beta_*}<q$, this implies that $\varkappa _{\alpha_*,\beta_*} \ge \varkappa _{\gamma,\beta_*}$. By the conditions $p_{\beta_*}<q\le p_\gamma$, we get
$$
\nu _{\alpha_*}\varkappa _{\alpha_*,\beta_*}^{1/p_\gamma‐1/p_{\alpha_*}} \stackrel{(\ref{kappa_ab})}{=} \nu _{\beta_*}\varkappa _{\alpha_*,\beta_*}^{1/p_\gamma‐1/p_{\beta_*}} \le 
\nu _{\beta_*}\varkappa _{\gamma,\beta_*}^{1/p_\gamma‐1/p_{\beta_*}} \stackrel{(\ref{kappa_ab})}{=} \nu_\gamma.
$$

This concludes the proof of (\ref{nu_a_g}). Hence,
$$
d_n(M, \, l_q^N) \stackrel{(\ref{nu_k_m})}{\gtrsim} d_n (\nu _{\alpha_*}k^{‐1/p_{\alpha_*}}V_k, \, l_q^N) \stackrel{(\ref{gl_q_l_2})}{\gtrsim} \nu_{\alpha_*}k^{1/q‐1/p_{\alpha_*}} \asymp \nu_{\alpha_*}\varkappa _{\alpha_*,\beta_*}^{1/q‐1/p_{\alpha_*}}.
$$

Let now $p_{\alpha_*}=q$. Then $\nu_{\alpha_*}\varkappa _{\alpha_*,\beta_*}^{1/q‐1/p_{\alpha_*}}=\nu_{\alpha_*}$. 

First we consider the case $\inf \{p_\alpha:\; \alpha \in A_1\}>q$. We take $\alpha_0\in A_1$ such that
\begin{align}
\label{a0_def}
\varkappa _{\alpha_0,\alpha_*} = \max _{\gamma\in A_1} \varkappa _{\gamma,\alpha_*}
\end{align}
(the maximum is attained, since $\{(\nu_\gamma, \, 1/p_\gamma):\; \gamma \in A_1\}$ is compact and we can apply (\ref{lim_kappa})).

We have
\begin{align}
\label{p_a_0_q}
p_{\alpha_0}>q. 
\end{align}

Let $k=[\varkappa_{\alpha_0,\alpha_*}]$. By (\ref{kappa_n}), we have $1\le k\le N$. We check that
\begin{align}
\label{nu_alpha_st_3case}
\nu_{\alpha_0} k^{‐1/p_{\alpha_0}}V_k \subset 2M.
\end{align}
It suffices to prove that 
\begin{align}
\label{nu_a_0}
\nu_{\alpha_0}\varkappa _{\alpha_0,\alpha_*}^{1/p_\gamma ‐ 1/p_{\alpha_0}} \le \nu_\gamma, \quad \gamma \in A.
\end{align}
Let $\gamma\in A_2'$. Then
$$
\nu_{\alpha_0}\varkappa _{\alpha_0,\alpha_*} ^{1/q ‐1/p_{\alpha_0}} = \nu_{\alpha_0}\varkappa _{\alpha_0,\alpha_*} ^{1/p_{\alpha_*} ‐1/p_{\alpha_0}} \stackrel{(\ref{kappa_ab})}{=} \nu _{\alpha_*} \stackrel{(\ref{nu_a_min})}{\le} \nu_{\alpha_0}\varkappa _{\alpha_0,\gamma} ^{1/q ‐1/p_{\alpha_0}};
$$
this together with (\ref{p_a_0_q}) yields that $\varkappa _{\alpha_0,\alpha_*} \le \varkappa _{\alpha_0,\gamma}$. Since $p_\gamma\le q<p_{\alpha_0}$, we have $$\nu_{\alpha_0}\varkappa _{\alpha_0,\alpha_*} ^{1/p_\gamma ‐1/p_{\alpha_0}} \le \nu_{\alpha_0}\varkappa _{\alpha_0,\gamma} ^{1/p_\gamma ‐1/p_{\alpha_0}} \stackrel{(\ref{kappa_ab})}{=}\nu_\gamma.$$

Let $\gamma\in A_1$. By (\ref{kappa_ab}), the inequality (\ref{nu_a_0}) is equivalent to $$\nu_{\alpha_*}\varkappa _{\alpha_0,\alpha_*} ^{1/p_\gamma ‐1/p_{\alpha_*}} \le \nu_{\alpha_*}\varkappa _{\gamma,\alpha_*} ^{1/p_\gamma ‐1/p_{\alpha_*}}.$$ Since $p_\gamma>q=p_{\alpha_*}$, it is equivalent to $\varkappa _{\alpha_0,\alpha_*}\ge \varkappa _{\gamma,\alpha_*}$. The last inequality follows from (\ref{a0_def}).

This completes the proof of (\ref{nu_a_0}). Hence,
\begin{align}
\label{dnm_lqn_nua}
d_n(M, \, l_q^N) \stackrel{(\ref{nu_alpha_st_3case})}{\gtrsim} \nu_{\alpha_0} k^{‐1/p_{\alpha_0}} d_n(V_k, \, l_q^N) \stackrel{(\ref{gl_q_l_2})}{\gtrsim} \nu_{\alpha_0} \varkappa_{\alpha_0,\alpha_*}^{1/q‐1/p_{\alpha_0}} \stackrel{(\ref{kappa_ab})}{=}\nu_{\alpha_*}\varkappa_{\alpha_0,\alpha_*}^{1/q‐1/p_{\alpha_*}}=\nu_{\alpha_*}.
\end{align}

Now let $\inf \{p_\alpha:\; \alpha\in A_1\}=q$. For $\delta>0$ we set 
\begin{align}
\label{a_delta}
A^\delta = A_2'\cup \{\alpha\in A_1:\; 1/q‐1/p_\alpha\ge \delta\}, 
\end{align}
$M_\delta =\cap _{\alpha\in A^\delta}\nu_\alpha B^N_{p_\alpha}$. By (\ref{dnm_lqn_nua}), we have
$$
d_n(M_\delta, \, l_q^N) \gtrsim \nu_{\alpha_*}.
$$
We prove that for small $\delta>0$ the inclusion $M_\delta\subset 2M$ holds. It suffices to take 
\begin{align}
\label{delta_log_n} 0<\delta \le \frac{\ln 2}{\ln N}
\end{align}
and to check that 
\begin{align}
\label{2a_a_a_del}
2\cap _{\alpha \in A\backslash A^\delta} \nu_\alpha B^N_{p_\alpha} \supset \nu_{\alpha_*}B^N_q.
\end{align}
Notice that $p_\alpha>q$ for $\alpha \in A\backslash A^\delta$. Hence (\ref{2a_a_a_del}) holds if $2\nu_\alpha\ge \nu_{\alpha_*}$ for all $\alpha\in A\backslash A^\delta$. We have $$\nu_{\alpha_*}/\nu_\alpha \stackrel{(\ref{kappa_ab})}{=}\varkappa_{\alpha,\alpha_*}^{1/p_{\alpha_*}‐1/p_\alpha} \stackrel{(\ref{kappa_n})}{\le} N^{1/q‐1/p_\alpha}\stackrel{(\ref{a_delta})}{\le} N^\delta\stackrel{(\ref{delta_log_n})}{\le} 2;$$
this implies (\ref{2a_a_a_del}) and completes the proof of the lower estimate of $d_n(M, \, l_q^N)$ for $p_{\alpha_*}=q$.

Let now $p_{\beta_*}=q$. Then $\nu_{\alpha_*}\varkappa _{\alpha_*,\beta_*}^{1/q‐1/p_{\alpha_*}} \stackrel{(\ref{kappa_ab})}{=} \nu _{\beta_*}\varkappa _{\alpha_*,\beta_*}^{1/q‐1/p_{\beta_*}} = \nu _{\beta_*}$. Hence, $d_n(M, \, l_q^N) \gtrsim \nu_{\beta_*}$.

{\bf Case 5.} We prove the upper estimate. Let $\alpha\in A_2'$. By Theorem \ref{glus}, $$d_n(M, \, l_q^N)\le d_n(\nu_\alpha B^N_{p_\alpha}, \, l_q^N) \underset{q}{\lesssim} \nu_\alpha (\min\{1, \, n^{‐1/2}N^{1/q}\})^{\frac{1/p_\alpha‐1/q}{1/2‐1/q}}.$$
Let $\alpha \in A_1'$, $\beta \in A_2'\cup A_3'$. Then $$M\subset \nu _\alpha B_{p_\alpha}^N \cap \nu_\beta B_{p_\beta}^N \stackrel{(\ref{kappa_ab}),(\ref{s_p_alpha})}{\subset} \nu_\alpha \varkappa _{\alpha,\beta} ^{1/q‐1/p_\alpha}B_q^N;$$
hence,
$$
d_n(M, \, l_q^N) \le \nu_\alpha \varkappa _{\alpha,\beta} ^{1/q‐1/p_\alpha}d_n(B_q^N, \, l_q^N) = \nu_\alpha \varkappa _{\alpha,\beta} ^{1/q‐1/p_\alpha}.
$$
Let $\alpha \in A_1'\cup A_2'$, $\beta \in A_3'$. Then $$M\subset \nu _\alpha B_{p_\alpha}^N \cap \nu_\beta B_{p_\beta}^N \stackrel{(\ref{kappa_ab}),(\ref{s_p_alpha})}{\subset} \nu_\alpha \varkappa _{\alpha,\beta} ^{1/2‐1/p_\alpha}B_2^N;$$
by Theorem \ref{glus},
$$
d_n(M, \, l_q^N) \le d_n(\nu_\alpha \varkappa _{\alpha,\beta} ^{1/2‐1/p_\alpha}B_2^N, \, l_q^N) \underset{q}{\lesssim} \nu_\alpha \varkappa _{\alpha,\beta} ^{1/2‐1/p_\alpha} \min \{1, \, n^{‐1/2}N^{1/q}\}.
$$

Now we prove the lower estimate. We denote $A_1=\{\alpha\in A:\; p_\alpha >q\}$, $A_2=\{\alpha\in A:\; 2<p_\alpha <q\}$, $A_3=\{\alpha\in A:\; p_\alpha <2\}$.
\begin{enumerate}
\item Let 
\begin{align}
\label{phi_a2}
\Phi(n, \, N, \, q, \, K)=\nu_{\alpha_*} (\min \{1, \, n^{‐1/2}N^{1/q}\})^{\frac{1/p_{\alpha_*}‐1/q}{1/2‐1/q}}, \quad \alpha_*\in A_2'; 
\end{align}
in addition, we suppose that $p_{\alpha_*}\ne q$ if $A_1\ne \varnothing$, and $p_{\alpha_*}\ne 2$ if $A_3\ne \varnothing$.

 We set 
 \begin{align}
 \label{l_max_nn}
 l = (\max\{1, \, n^{1/2}N^{‐1/q}\})^{\frac{1}{1/2‐1/q}}, \;\; k = \lceil l \rceil.
 \end{align}
Notice that $1\le l\le N$ and $1\le k\le N$. We prove that 
\begin{align}
\label{nu_a_k_2m}
\nu_{\alpha_*} k^{‐1/p_{\alpha_*}}V_k \subset 2M. 
\end{align}
It suffices to check that 
$$
\nu_{\alpha_*} l^{1/p_\gamma‐1/p_{\alpha_*}} \le \nu _\gamma, \quad \gamma \in A.
$$
If $\gamma\in A_2'$, it follows from the inequality
$$
\nu_{\alpha_*}(\min \{1, \, n^{‐1/2}N^{1/q}\}) ^{\frac{1/p_{\alpha_*}‐1/q}{1/2‐1/q}} \stackrel{(\ref{5}), (\ref{7}),(\ref{phi_a2})}{\le} \nu_{\gamma}(\min \{1, \, n^{‐1/2}N^{1/q}\}) ^{\frac{1/p_{\gamma}‐1/q}{1/2‐1/q}}.
$$
Let $\gamma \in A_1$ (hence, $A_1\ne \varnothing$; by our assumption, we have $p_{\alpha_*}<q$). Then
$$
\nu_{\alpha_*}(\min \{1, \, n^{‐1/2}N^{1/q}\}) ^{\frac{1/p_{\alpha_*}‐1/q}{1/2‐1/q}} \stackrel{(\ref{5}), (\ref{6}), (\ref{phi_a2})}{\le} \nu_{\gamma} \varkappa _{\gamma,\alpha_*}^{1/q‐1/p_\gamma}\stackrel{(\ref{kappa_ab})}{=} \nu_{\alpha_*}\varkappa _{\gamma,\alpha_*}^{1/q‐1/p_{\alpha_*}}.
$$
Hence, $l^{1/q‐1/p_{\alpha_*}}\le \varkappa _{\gamma,\alpha_*}^{1/q‐1/p_{\alpha_*}}$, and $l\ge \varkappa _{\gamma,\alpha_*}$. Since $p_\gamma >q>p_{\alpha_*}$, we have
$$
\nu_{\alpha_*}l^{1/p_\gamma‐1/p_{\alpha_*}}\le \nu_{\alpha_*}\varkappa _{\gamma,\alpha_*}^{1/p_\gamma‐1/p_{\alpha_*}} \stackrel{(\ref{kappa_ab})}{=} \nu_\gamma.
$$
Let $\gamma \in A_3$ (hence, $A_3\ne \varnothing$; by our assumption, we have $p_{\alpha_*}>2$). Then
$$
\nu_{\alpha_*}(\min \{1, \, n^{‐1/2}N^{1/q}\}) ^{\frac{1/p_{\alpha_*}‐1/q}{1/2‐1/q}} \stackrel{(\ref{5}), (\ref{8}),(\ref{phi_a2})}{\le} \nu_{\alpha_*} \varkappa _{\alpha_*,\gamma}^{1/2‐1/p_{\alpha_*}}\min \{1, \, n^{‐1/2}N^{1/q}\}.
$$
Therefore, $l^{1/2‐1/p_{\alpha_*}}\le \varkappa _{\alpha_*,\gamma}^{1/2‐1/p_{\alpha_*}}$. Hence, $l\le \varkappa _{\alpha_*,\gamma}$; since $p_\gamma <2<p_{\alpha_*}$, we have
$$
\nu_{\alpha_*}l^{1/p_\gamma‐1/p_{\alpha_*}}\le 
\nu_{\alpha_*}\varkappa _{\alpha_*,\gamma}^{1/p_\gamma‐1/p_{\alpha_*}}\stackrel{(\ref{kappa_ab})}{=} \nu_\gamma.
$$

This concludes the proof of (\ref{nu_a_k_2m}). Notice that $n\stackrel{(\ref{l_max_nn})}{\le} N^{2/q}k^{1‐2/q}$. Hence,
$$d_n(M, \, l_q^N) \stackrel{(\ref{nu_a_k_2m})}{\gtrsim} d_n (\nu_{\alpha_*} k^{‐1/p_{\alpha_*}}V_k, \, l_q^N) \stackrel{(\ref{kq1})}{\underset{q}{\gtrsim}} \nu_{\alpha_*}k^{1/q‐1/p_{\alpha_*}} \stackrel{(\ref{l_max_nn})}{\gtrsim} $$$$ \gtrsim \nu _{\alpha_*}(\min \{1, \, n^{‐1/2}N^{1/q}\})^{\frac{1/p_{\alpha_*}‐1/q}{1/2‐1/q}}.$$

\item Let 
\begin{align}
\label{phi_a12} \Phi(n, \, N, \, q, \, K) = \nu_{\alpha_*}, \quad p_{\alpha_*}=q, \quad A_1\ne \varnothing.
\end{align}

First we suppose that 
\begin{align}
\label{sup_pa_l_q}
\sup \{p_\alpha:\; \alpha \in A_2\cup A_3'\}< q.
\end{align}
 We take $\alpha_0\in A_2 \cup A_3'$ such that 
\begin{align}
\label{kappa_a0ast_min}
\varkappa _{\alpha_*,\alpha_0} = \min _{\gamma\in A_2\cup A_3'} \varkappa _{\alpha_*,\gamma}
\end{align}
 (by (\ref{sup_pa_l_q}), the minimun is attained and $p_{\alpha_0}<q$). Let $k=\lceil \varkappa_{\alpha_*,\alpha_0}\rceil$. From (\ref{kappa_n}) it follows that $1\le k\le N$. We prove that 
\begin{align}
\label{nu_a0_kp0}
\nu_{\alpha_0} k^{‐1/p_{\alpha_0}}V_k \subset 2M. 
\end{align}
It suffices to check the inequality 
\begin{align}
\label{nu_a0_kappaast}
\nu_{\alpha_0} \varkappa _{\alpha_*,\alpha_0} ^{1/p_\gamma ‐ 1/p_{\alpha_0}}\le \nu_\gamma, \quad \gamma\in A.
\end{align}
Let $\gamma \in A_2\cup A_3'$. From (\ref{kappa_a0ast_min}) it follows that $\varkappa _{\alpha_*,\alpha_0} \le \varkappa _{\alpha_*,\gamma}$. Since $p_{\alpha_*}=q > p_\gamma$, we get $\nu_{\alpha_*}\varkappa _{\alpha_*,\alpha_0} ^{1/p_\gamma ‐ 1/p_{\alpha_*}}\le \nu_{\alpha_*} \varkappa _{\alpha_*,\gamma} ^{1/p_\gamma‐1/p_{\alpha_*}}$. By (\ref{kappa_ab}), this inequality is equivalent to (\ref{nu_a0_kappaast}).

Let $\gamma\in A_1'$. By (\ref{kappa_ab}), the inequality (\ref{nu_a0_kappaast}) is equivalent to $\nu_{\alpha_0} \varkappa _{\alpha_*,\alpha_0}^{1/p_\gamma‐1/p_{\alpha_0}}\le \nu_{\alpha_0} \varkappa _{\gamma,\alpha_0}^{1/p_\gamma‐1/p_{\alpha_0}}$. Since $p_{\alpha_0}<p_\gamma$, it is equivalent to 
\begin{align}
\label{kappa_ast_a0_kappa_g_a0}
\varkappa _{\alpha_*,\alpha_0}\ge \varkappa _{\gamma,\alpha_0}. 
\end{align}
Let us prove this inequality. If $p_\gamma=q$, then $\gamma=\alpha_*$. Let $p_\gamma>q$. We have $$\nu_{\alpha_*} \stackrel{(\ref{5}),(\ref{6}),(\ref{phi_a12})}{\le} \nu _{\gamma} \varkappa _{\gamma,\alpha_0}^{1/q‐1/p_\gamma};$$ by (\ref{kappa_ab}), it is equivalent to $\nu_{\alpha_0}\varkappa _{\alpha_*,\alpha_0}^{1/p_{\alpha_*}‐1/p_{\alpha_0}} \le \nu_{\alpha_0} \varkappa _{\gamma,\alpha_0}^{1/q‐1/p_{\alpha_0}} = \nu_{\alpha_0} \varkappa _{\gamma,\alpha_0}^{1/p_{\alpha_*}‐1/p_{\alpha_0}}$. Since $p_{\alpha_0}<q =p_{\alpha_*}$, we get (\ref{kappa_ast_a0_kappa_g_a0}).

This completes the proof of (\ref{nu_a0_kp0}).

Now we check that 
\begin{align}
\label{nn2qkaa0}
n\le N^{2/q}\varkappa _{\alpha_*,\alpha_0}^{1‐2/q}.
\end{align}
It suffices to consider $n\ge N^{2/q}$.

If $\alpha_0\in A_2$, (\ref{nn2qkaa0}) follows from
$$
\nu_{\alpha_*} \stackrel{(\ref{5}),(\ref{7}),(\ref{phi_a12})}{\le} \nu_{\alpha_0}(n^{‐1/2}N^{1/q})^{\frac{1/p_{\alpha_0}‐1/q}{1/2‐1/q}} \stackrel{(\ref{kappa_ab})}{=} \nu_{\alpha_*} \varkappa_{\alpha_*,\alpha_0} ^{1/p_{\alpha_0} ‐1/p_{\alpha_*}}(n^{‐1/2}N^{1/q})^{\frac{1/p_{\alpha_0}‐1/q}{1/2‐1/q}}.
$$
If $\alpha_0\in A_3'$, (\ref{nn2qkaa0}) follows from
$$
\nu_{\alpha_*} \stackrel{(\ref{5}), (\ref{8}), (\ref{phi_a12})}{\le} \nu_{\alpha_*} \varkappa _{\alpha_*,\alpha_0}^{1/2‐1/p_{\alpha_*}}n^{‐1/2}N^{1/q}.
$$
Hence,
\begin{align}
\label{dnmlqnnua}
d_n(M, \, l_q^N) \stackrel{(\ref{kq1}), (\ref{nu_a0_kp0}), (\ref{nn2qkaa0})}{\underset{q}{\gtrsim}} \nu _{\alpha_0}k^{1/q‐1/p_{\alpha_0}} \asymp \nu_{\alpha_0} \varkappa 
_{\alpha_*,\alpha_0}^{1/p_{\alpha_*}‐1/p_{\alpha_0}} \stackrel{(\ref{kappa_ab})}{=} \nu_{\alpha_*}.
\end{align}

Let now 
$$
\sup \{p_\alpha:\; \alpha \in A_2\cup A_3'\}= q.
$$
We set 
\begin{align}
\label{a_del_22}
A^\delta =A_1'\cup \{\alpha\in A_2\cup A_3':\; 1/p_\alpha‐1/q\ge \delta\},
\end{align}
$M_\delta = \cap _{\alpha \in A^\delta} \nu_\alpha B^N _{p_\alpha}$. By (\ref{dnmlqnnua}),
$$
d_n(M_\delta, \, l_q^N) \underset{q}{\gtrsim} \nu_{\alpha_*}.
$$
We prove that if 
\begin{align}
\label{delta_2N_2}
0<\delta \le \frac{\ln 2}{\ln N}, 
\end{align}
then $M_\delta \subset 2M$. It suffices to check that
\begin{align}
\label{2capaaadelbnq}
2\cap _{\alpha\in A\backslash A^\delta} \nu_\alpha B^N_{p_\alpha} \supset \nu_{\alpha_*}B^N_q.
\end{align}
Since $p_\alpha<q$ for all $\alpha\in A\backslash A^\delta$, (\ref{2capaaadelbnq}) holds if 
\begin{align}
\label{111111}
\nu_{\alpha_*}\le 2\nu_\alpha \cdot N^{1/q‐1/p_\alpha}, \quad \alpha \in A\backslash A^\delta. 
\end{align}
We have $\nu _{\alpha} \stackrel{(\ref{kappa_ab})}{=} \nu_{\alpha_*}\varkappa _{\alpha_*,\alpha}^{1/p_\alpha ‐1/p_{\alpha_*}} \stackrel{(\ref{kappa_n}), (\ref{a_del_22})}{\ge} \nu_{\alpha_*}$. This together with (\ref{a_del_22}) and (\ref{delta_2N_2}) yields (\ref{111111}).

\item Let 
\begin{align}
\label{phi_min_3} \Phi(n, \, N, \, q, \, K) = \nu_{\alpha_*} n^{‐1/2} N^{1/q}, \quad p_{\alpha_*}=2, \quad A_3\ne \varnothing,
\end{align}
and let $n\ge N^{2/q}$. 

First we consider the case
\begin{align}
\label{inf_pa_g_2}
\inf \{p_\alpha :\; \alpha \in A_1'\cup A_2\}>2. 
\end{align}
We take $\alpha_0\in A_1'\cup A_2$ such that 
\begin{align}
\label{kappa_a0a_max_kga}
\varkappa _{\alpha_0,\alpha_*} =\max _{\gamma\in A_1'\cup A_2} \varkappa _{\gamma,\alpha_*}
\end{align}
(by (\ref{inf_pa_g_2}), the maximum is attained and $p_{\alpha_0}>2$). Let $k=\lfloor \varkappa _{\alpha_0,\alpha_*}\rfloor$. 

We check that 
\begin{align}
\label{2nua0kpa02m}
\nu_{\alpha_0}k^{‐1/p_{\alpha_0}}V_k \subset 2M.
\end{align} 
It suffices to prove the inequality 
\begin{align}
\label{nua0vkppa0astnugam}
\nu_{\alpha_0}\varkappa _{\alpha_0,\alpha_*} ^{1/p_\gamma‐1/p_{\alpha_0}} \le \nu _\gamma, \quad \gamma \in A.
\end{align} 

If $\gamma\in A_3'$, then 
$$
\nu_{\alpha_0}\varkappa _{\alpha_0,\gamma}^{1/2‐1/p_{\alpha_0}}n^{‐1/2}N^{1/q}\stackrel{(\ref{5}),(\ref{8}),(\ref{phi_min_3})}{\ge} \nu_{\alpha_*} n^{‐1/2} N^{1/q} \stackrel{(\ref{kappa_ab})}{=}$$$$= \nu_{\alpha_0}\varkappa _{\alpha_0,\alpha_*}^{1/p_{\alpha_*}‐1/p_{\alpha_0}}n^{‐1/2}N^{1/q} = \nu_{\alpha_0}\varkappa _{\alpha_0,\alpha_*}^{1/2‐1/p_{\alpha_0}}n^{‐1/2}N^{1/q};
$$
this together with the condition $p_{\alpha_0}>2$ yields that $\varkappa _{\alpha_0,\gamma}\ge \varkappa _{\alpha_0,\alpha_*}$. Since $p_\gamma\le 2<p_{\alpha_0}$, we have
$$
\nu_{\alpha_0}\varkappa _{\alpha_0,\alpha_*} ^{1/p_\gamma‐1/p_{\alpha_0}} \le \nu_{\alpha_0}\varkappa _{\alpha_0,\gamma} ^{1/p_\gamma‐1/p_{\alpha_0}} \stackrel{(\ref{kappa_ab})}{=} \nu_\gamma.
$$
If $\gamma\in A_1'\cup A_2$, then by (\ref{kappa_ab})  inequality (\ref{nua0vkppa0astnugam}) is equivalent to $\nu_{\alpha_*}\varkappa _{\alpha_0,\alpha_*}^{1/p_\gamma‐1/p_{\alpha_*}} \le \nu_{\alpha_*}\varkappa _{\gamma,\alpha_*}^{1/p_\gamma‐1/p_{\alpha_*}}$; since $p_\gamma>p_{\alpha_*}$, it is equivalent to $\varkappa _{\alpha_0,\alpha_*} \ge \varkappa _{\gamma,\alpha_*}$; the last inequality holds by (\ref{kappa_a0a_max_kga}).

This completes the proof of (\ref{2nua0kpa02m}).

Let us prove that 
\begin{align}
\label{ngan2qka0a}
n\ge N^{2/q}\varkappa _{\alpha_0,\alpha_*}^{1‐2/q}. 
\end{align}
For $\alpha_0\in A_1'$, it follows from
$$
\nu_{\alpha_*}n^{‐1/2}N^{1/q} \stackrel{(\ref{5}), (\ref{6}), (\ref{phi_min_3})}{\le} \nu_{\alpha_0} \varkappa _{\alpha_0,\alpha_*}^{1/q‐1/p_{\alpha_0}} \stackrel{(\ref{kappa_ab})}{=} \nu_{\alpha_*} \varkappa _{\alpha_0,\alpha_*}^{1/q‐1/p_{\alpha_*}}=\nu_{\alpha_*} \varkappa _{\alpha_0,\alpha_*}^{1/q‐1/2},
$$
for $\alpha_0\in A_2$, it follows from
$$
\nu _{\alpha_*}n^{‐1/2}N^{1/q} \stackrel{(\ref{5}), (\ref{7}), (\ref{phi_min_3})}{\le} \nu_{\alpha_0} (n^{‐1/2} N^{1/q})^{\frac{1/p_{\alpha_0}‐1/q}{1/2‐1/q}} \stackrel{(\ref{kappa_ab})}{=} \nu_{\alpha_*} \varkappa_{\alpha_0,\alpha_*}^{1/p_{\alpha_0}‐1/p_{\alpha_*}} (n^{‐1/2} N^{1/q})^{\frac{1/p_{\alpha_0}‐1/q}{1/2‐1/q}}.
$$

Hence, we get $$d_n(M, \, l_q^N)\stackrel{(\ref{2nua0kpa02m})}{\gtrsim}d_n(\nu_{\alpha_0}k^{‐1/p_{\alpha_0}}V_k, \, l_q^N) \stackrel{(\ref{kq0}), (\ref{ngan2qka0a})}{\underset{q}{\gtrsim}} $$$$\gtrsim\nu_{\alpha_0}\varkappa _{\alpha_0,\alpha_*}^{1/2‐1/p_{\alpha_0}}n^{‐1/2}N^{1/q} \stackrel{(\ref{kappa_ab})}{=} \nu_{\alpha_*}n^{‐1/2}N^{1/q}.$$

If $\inf \{p_\alpha :\; \alpha \in A_1'\cup A_2\}=2$, we argue as in the previous cases and obtain the desired estimate.

\item Let 
\begin{align}
\label{nu_a_phi_gq}
\Phi(n, \, N, \, q, \, K)=\nu_{\alpha_*}\varkappa _{\alpha_*,\beta_*}^{1/q‐1/p_{\alpha_*}}, \quad \alpha_* \in A_1', \;\; \beta_* \in A_2'\cup A_3'.
\end{align}
Notice that if $p_{\alpha_*}=q$, the right‐hand side of (\ref{nu_a_phi_gq}) is equal to $\nu_{\alpha_*}$; if $p_{\beta_*}=q$, we get by (\ref{kappa_ab})  $\nu _{\alpha_*}\varkappa _{\alpha_*,\beta_*}^{1/q‐1/p_{\alpha_*}}=\nu_{\beta_*}$; thus, we have arrived at the case already considered. Further we suppose that $p_{\alpha_*}> q$, $p_{\beta_*}<q$.

Let $k =\lceil \varkappa _{\alpha_*,\beta_*}\rceil$.

We prove that 
\begin{align}
\label{nual1pa2p4}
\nu_{\alpha_*} k^{‐1/p_{\alpha_*}}V_k \subset 2M.
\end{align}
It suffices to check that
\begin{align}
\label{nual1pa2p41111}
\nu_{\alpha_*} \varkappa _{\alpha_*,\beta_*}^{1/p_\gamma‐1/p_{\alpha_*}} \le \nu_\gamma, \quad \gamma \in A.
\end{align}
Let $\gamma\in A_1'$. We have $\nu_{\alpha_*} \varkappa _{\alpha_*,\beta_*}^{1/q‐1/p_{\alpha_*}} \stackrel{(\ref{5}), (\ref{6}), (\ref{nu_a_phi_gq})}{\le} \nu_\gamma \varkappa _{\gamma,\beta_*}^{1/q ‐1/p_\gamma}$; by (\ref{kappa_ab}), it is equivalint to $\nu_{\beta_*} \varkappa _{\alpha_*,\beta_*}^{1/q‐1/p_{\beta_*}} \le \nu_{\beta_*} \varkappa _{\gamma,\beta_*}^{1/q ‐1/p_{\beta_*}}$; since $p_{\beta_*}<q$, we have $\varkappa _{\alpha_*,\beta_*} \ge \varkappa _{\gamma, \beta_*}$. This together with $p_\gamma \ge q>p_{\beta_*}$ implies that
$$
\nu_{\alpha_*} \varkappa _{\alpha_*,\beta_*}^{1/p_\gamma‐1/p_{\alpha_*}} \stackrel{(\ref{kappa_ab})}{=} \nu_{\beta_*} \varkappa _{\alpha_*,\beta_*}^{1/p_\gamma‐1/p_{\beta_*}} \le \nu_{\beta_*} \varkappa _{\gamma,\beta_*}^{1/p_\gamma‐1/p_{\beta_*}} \stackrel{(\ref{kappa_ab})}{=} \nu_\gamma.
$$
Let $\gamma\in A_2'\cup A_3'$. We have
$$
\nu_{\alpha_*} \varkappa _{\alpha_*,\beta_*}^{1/q‐1/p_{\alpha_*}} \stackrel{(\ref{5}), (\ref{6}), (\ref{nu_a_phi_gq})}{\le} \nu_{\alpha_*} \varkappa _{\alpha_*,\gamma}^{1/q‐1/p_{\alpha_*}};
$$
hence, $\varkappa _{\alpha_*,\beta_*} \le \varkappa _{\alpha_*,\gamma}$ by the condition $p_{\alpha_*}>q$. Since $p_\gamma\le q< p_{\alpha_*}$, we have
$$
\nu_{\alpha_*} \varkappa _{\alpha_*,\beta_*}^{1/p_\gamma‐1/p_{\alpha_*}} \le \nu_{\alpha_*} \varkappa _{\alpha_*,\gamma}^{1/p_\gamma‐1/p_{\alpha_*}} \stackrel{(\ref{kappa_ab})}{=} \nu_\gamma.
$$

This concludes the proof of (\ref{nual1pa2p41111}).

Now we check the inequality 
\begin{align}
\label{nlen2qkab}
n\le N^{2/q}\varkappa _{\alpha_*,\beta_*}^{1‐2/q}. 
\end{align}
It suffices to consider $n\ge N^{2/q}$. 

If $\beta_*\in A_2$, then (\ref{nlen2qkab}) follows from
$$
\nu_{\beta_*}\varkappa _{\alpha_*,\beta_*}^{1/q‐1/p_{\beta_*}} \stackrel{ (\ref{kappa_ab})}{=}\nu_{\alpha_*}\varkappa _{\alpha_*,\beta_*}^{1/q‐1/p_{\alpha_*}} \stackrel{(\ref{5}), (\ref{7}),(\ref{nu_a_phi_gq})}{\le} \nu_{\beta_*}(n^{‐1/2}N^{1/q}) ^{\frac{1/p_{\beta_*}‐1/q}{1/2‐1/q}}.
$$
If $\beta_*\in A_3'$, then (\ref{nlen2qkab}) follows from
$$
\nu_{\alpha_*}\varkappa _{\alpha_*,\beta_*}^{1/q‐1/p_{\alpha_*}} \stackrel{(\ref{5}), (\ref{8}), (\ref{nu_a_phi_gq})}{\le} \nu_{\alpha_*}\varkappa _{\alpha_*,\beta_*} ^{1/2‐1/p_{\alpha_*}}n^{‐1/2}N^{1/q}.
$$

Hence,
$$
d_n(M, \, l_q^N)\stackrel{(\ref{nual1pa2p4})}{\gtrsim} d_n(\nu_{\alpha_*} k^{‐1/p_{\alpha_*}}V_k, \, l_q^N) \stackrel{(\ref{kq1}),(\ref{nlen2qkab})}{\underset{q}{\gtrsim}} $$$$\gtrsim\nu_{\alpha_*} k^{1/q‐1/p_{\alpha_*}} \gtrsim \nu_{\alpha_*} \varkappa _{\alpha_*,\beta_*} ^{1/q‐1/p_{\alpha_*}}.
$$

\item Let 
\begin{align}
\label{phi_last}
\Phi(n, \, N, \, q, \, K)=\nu_{\alpha_*} \varkappa _{\alpha_*,\beta_*}^{1/2‐1/p_{\alpha_*}}n^{‐1/2}N^{1/q},\quad \alpha_*\in A_1'\cup A_2', \quad \beta_*\in A_3', 
\end{align}
and let $n\ge N^{2/q}$. Notice that if $p_{\alpha_*}=2$, the right‐hand side of (\ref{phi_last}) equals to $\nu_{\alpha_*}n^{‐1/2}N^{1/q}$; if $p_{\beta_*}=2$, we have $\nu_{\alpha_*} \varkappa _{\alpha_*,\beta_*}^{1/2‐1/p_{\alpha_*}}n^{‐1/2}N^{1/q} \stackrel{(\ref{kappa_ab})}{=}\nu_{\beta_*}n^{‐1/2}N^{1/q}$. Thus, we have arrived at the case already considered.

Further we suppose that $p_{\alpha_*}>2$, $p_{\beta_*}<2$.

Let $k =\lfloor \varkappa _{\alpha_*,\beta_*}\rfloor$. 

We show that 
\begin{align}
\label{nu_a_k_vk_2m}
\nu_{\alpha_*}k^{‐1/p_{\alpha_*}}V_k \subset 2M. 
\end{align}
It suffices to check that
\begin{align}
\label{nuastvarkaabb}
\nu_{\alpha_*} \varkappa _{\alpha_*,\beta_*} ^{1/p_\gamma‐1/p_{\alpha_*}} \le \nu_\gamma, \quad \gamma\in A.
\end{align}
If $\gamma\in A_1'\cup A_2'$, we have
$$
\nu_{\alpha_*}\varkappa _{\alpha_*,\beta_*}^{1/2‐1/p_{\alpha_*}}n^{‐1/2}N^{1/q} \stackrel{(\ref{5}), (\ref{8}),(\ref{phi_last})}{\le} \nu_{\gamma}\varkappa _{\gamma,\beta_*}^{1/2‐1/p_{\gamma}}n^{‐1/2}N^{1/q};
$$
this together with (\ref{kappa_ab}) yields that $\nu_{\beta_*}\varkappa _{\alpha_*,\beta_*}^{1/2‐1/p_{\beta_*}}\le \nu_{\beta_*}\varkappa _{\gamma,\beta_*}^{1/2‐1/p_{\beta_*}}$. Since $p_{\beta_*}<2$, it implies that $\varkappa _{\alpha_*,\beta_*} \ge \varkappa _{\gamma,\beta_*}$. Hence, by $p_\gamma \ge 2>p_{\beta_*}$ we get
$$
\nu_{\alpha_*} \varkappa _{\alpha_*,\beta_*} ^{1/p_\gamma‐1/p_{\alpha_*}} \stackrel{(\ref{kappa_ab})}{=} \nu_{\beta_*} \varkappa _{\alpha_*,\beta_*} ^{1/p_\gamma‐1/p_{\beta_*}} \le \nu_{\beta_*} \varkappa _{\gamma,\beta_*} ^{1/p_\gamma‐1/p_{\beta_*}} \stackrel{(\ref{kappa_ab})}{=} \nu_\gamma.
$$
If $\gamma \in A_3'$, we have $$\nu_{\alpha_*}\varkappa _{\alpha_*,\beta_*}^{1/2‐1/p_{\alpha_*}}n^{‐1/2}N^{1/q} \stackrel{(\ref{5}), (\ref{8}),(\ref{phi_last})}{\le} \nu_{\alpha_*} \varkappa _{\alpha_*,\gamma}^{1/2‐1/p_{\alpha_*}}n^{‐1/2}N^{1/q}.$$ Since $p_{\alpha_*}>2$, it yields that $\varkappa _{\alpha_*,\beta_*}\le \varkappa _{\alpha_*,\gamma}$. Recall that $p_{\alpha_*}>2\ge p_\gamma$. Hence,
$$
\nu_{\alpha_*} \varkappa _{\alpha_*,\beta_*} ^{1/p_\gamma‐1/p_{\alpha_*}} \le \nu_{\alpha_*} \varkappa _{\alpha_*,\gamma} ^{1/p_\gamma‐1/p_{\alpha_*}} \stackrel{(\ref{kappa_ab})}{=} \nu_\gamma.
$$

This completes the proof of (\ref{nuastvarkaabb}).

Now we show that 
\begin{align}
\label{ngevrkab12qn2q}
n\ge \varkappa_{\alpha_*,\beta_*} ^{1‐2/q}N^{2/q}.
\end{align}

If $\alpha_*\in A_1'$, then (\ref{ngevrkab12qn2q}) follows from $$\nu_{\alpha_*}\varkappa _{\alpha_*,\beta_*}^{1/2‐1/p_{\alpha_*}}n^{‐1/2}N^{1/q} \stackrel{(\ref{5}), (\ref{6}),(\ref{phi_last})}{\le} \nu_{\alpha_*} \varkappa _{\alpha_*,\beta_*} ^{1/q‐1/p_{\alpha_*}}.$$ If $\alpha_*\in A_2$, then (\ref{ngevrkab12qn2q}) follows from
$$
\nu_{\alpha_*}\varkappa _{\alpha_*,\beta_*}^{1/2‐1/p_{\alpha_*}}n^{‐1/2}N^{1/q} \stackrel{(\ref{5}), (\ref{7}), (\ref{phi_last})}{\le} \nu_{\alpha_*}(n^{‐1/2}N^{1/q})^{\frac{1/p_{\alpha_*}‐1/q}{1/2‐1/q}}.
$$ 

Hence, $$d_n(M, \, l_q^N) \stackrel{(\ref{nu_a_k_vk_2m})}{\gtrsim} d_n(\nu_{\alpha_*}k^{‐1/p_{\alpha_*}}V_k, \, l_q^N) \stackrel{(\ref{kq0}),(\ref{ngevrkab12qn2q})}{\underset{q}{\gtrsim}}$$$$\gtrsim \nu_{\alpha_*} k^{1/2‐1/p_{\alpha_*}}n^{‐1/2}N^{1/q} \asymp \nu_{\alpha_*} \varkappa _{\alpha_*,\beta_*}^{1/2‐1/p_{\alpha_*}}n^{‐1/2}N^{1/q}.$$

\item Let $n\le N^{2/q}$, $\Phi(n, \, N, \, q, \, K) = \nu_{\alpha_*}\varkappa _{\alpha_*,\beta_*}^{1/2‐1/p_{\alpha_*}}$, $\alpha_*\in A_1'\cup A_2$, $\beta_*\in A_3$. If $\alpha_*\in A_1'$, it follows from (\ref{5}), (\ref{6}) that $\nu_{\alpha_*}\varkappa _{\alpha_*,\beta_*}^{1/2‐1/p_{\alpha_*}} \le \nu_{\alpha_*}\varkappa _{\alpha_*,\beta_*}^{1/q‐1/p_{\alpha_*}}$; therefore, $\varkappa _{\alpha_*,\beta_*}\le 1$. By (\ref{kappa_n}), we get $\varkappa _{\alpha_*,\beta_*}=1$ and $\Phi(n, \, N, \, q, \, K) = \nu_{\alpha_*}\varkappa _{\alpha_*,\beta_*}^{1/q‐1/p_{\alpha_*}}$. We have arrived at the case already considered. If $\alpha_*\in A_2$, it follows from (\ref{5}), (\ref{7}) that $\nu_{\alpha_*}\varkappa _{\alpha_*,\beta_*}^{1/2‐1/p_{\alpha_*}} \le \nu_{\alpha_*}$. Hence, $\varkappa _{\alpha_*,\beta_*}\stackrel{(\ref{kappa_n})}{=}1$ and $\Phi(n, \, N, \, q, \, K) = \nu_{\alpha_*}$. We have arrived at the case already considered.

\item Let $n\le N^{2/q}$, $p_{\alpha_*}= 2$, $\Phi(n, \, N, \, q, \, K) = \nu_{\alpha_*}$. We prove that $\nu_{\alpha_*}\le \nu_\alpha$ for all $\alpha \in A$. If $\alpha\in A_2$, it follows from (\ref{5}), (\ref{7}). If $\alpha\in A_3'$, then $\nu_\alpha \stackrel{(\ref{kappa_ab})}{=}\nu_{\alpha_*}\varkappa _{\alpha_*,\alpha}^{1/p_\alpha‐1/p_{\alpha_*}}\stackrel{(\ref{kappa_n})}{\ge} \nu_{\alpha_*}$. If $\alpha\in A_1'$, then $\nu_{\alpha_*}\le \nu_\alpha \varkappa_{\alpha,\alpha_*}^{1/q‐1/p_{\alpha}}$ by (\ref{5}), (\ref{6}). Applying (\ref{kappa_ab}), we get $\nu_{\alpha_*}\le \nu_{\alpha_*} \varkappa_{\alpha,\alpha_*}^{1/q‐1/p_{\alpha_*}}$, or $\varkappa_{\alpha,\alpha_*}\le 1$. By (\ref{kappa_n}), we have $\varkappa_{\alpha,\alpha_*}=1$ and $\nu_{\alpha_*}= \nu_\alpha$.  

Hence, $M \supset \nu_{\alpha_*}B_1^N$, and by Theorem \ref{glus} $$d_n(M, \, l_q^N) \ge \nu_{\alpha_*}d_n(B_1^N, \, l_q^N) \underset{q}{\gtrsim} \nu_{\alpha_*}.$$

\end{enumerate}
This completes the proof.
\end{proof}

\renewcommand{\proofname}{\bf Proof of Theorem \ref{main_infty}}

\begin{proof}
The upper estimate follows from Theorem \ref{garn_glusk_teor}.

Let us prove the lower estimate. As in Theorem \ref{main}, we may assume that the set $K$ is compact and the minimum in the right‐hand side of (\ref{d_n_infty_trm}) is attained at some $\alpha_*\in A$. Then 
\begin{align}
\label{infty_ineq} \nu_{\alpha_*} \left\lceil n \log ^{‐1} \frac {2N}{n} \right \rceil ^{‐1/p_{\alpha_*}} \le \nu_{\alpha} \left\lceil n \log ^{‐1} \frac {2N}{n} \right \rceil ^{‐1/p_{\alpha}}, \quad \alpha\in A.
\end{align}
Let $C>0$ be the constant from Theorem \ref{t_kmr}, $k_* = \min \left\{\left\lceil C n \log ^{‐1} \frac {2N}{n} \right \rceil, \, N \right\}$. Notice that if $k_*=N$, then $n\asymp N$, $d_n(V_{k_*}, \, l_\infty^N) = d_n(B_\infty ^N, \, l_\infty^N) = 1$.

We have $$d_n (\nu_{\alpha_*}k_*^{‐1/p_{\alpha_*}}V_{k_*}, \, l_\infty^N) \ge \nu_{\alpha_*}k_*^{‐1/p_{\alpha_*}}/2 \gtrsim \nu_{\alpha_*} \min \left\{1, \, n^{‐1}\log \frac {2N}{n}\right\}^{1/p_{\alpha_*}}.$$
By (\ref{infty_ineq}), there is an absolute constant $C_1>0$ such that $\nu_{\alpha_*}k_*^{‐1/p_{\alpha_*}}V_{k_*} \subset C_1M$. This concludes the proof.
\end{proof}

\section{The case of arbitrary $\nu_\alpha$.}

Notice that if $\inf _{\alpha\in A} \nu_\alpha=0$, then $\cap _{\alpha\in A} \nu_\alpha B^N_{p_\alpha}=\{0\}$. Therefore, further we assume that $\inf _{\alpha\in A} \nu_\alpha>0$.

We set $Z=\{1/p_\alpha:\; \alpha\in A\} \subset [0, \, 1]$; for $z=1/p_\alpha$ we set $\nu(z)=\nu_\alpha$. Then $M = \cap _{z\in Z} \nu(z) B^N_{1/z}$, and (\ref{kappa_n}) is equivalent to
$$
1\le \frac{\nu(z)}{\nu(w)} \le N^{z‐w}, \quad z, \, w\in Z, \quad z\ge w.
$$

Let $\nu:Z \rightarrow (0, \, \infty)$ be an arbitrary function, and let $\inf _{z\in Z} \nu(z)>0$. We construct the function $\nu_{**}:Z \rightarrow (0, \, \infty)$ such that
\begin{align}
\label{mon_lip} 1\le \frac{\nu_{**}(z)}{\nu_{**}(w)} \le N^{z‐w}, \quad z, \, w\in Z, \quad z\ge w,
\end{align}
\begin{align}
\label{m_cap_zz}
M=\cap _{z\in Z} \nu_{**}(z) B^N_{1/z}.
\end{align}

Let 
$$
\nu_*(z) = \inf \{\nu(w): w\ge z\}, \quad z\in Z.
$$
Then $\nu_*(z)\le \nu(z)$ and the function $\nu_*$ is non‐decreasing. We show that
\begin{align}
\label{m_nu_st} M = \cap _{z\in Z} \nu_*(z) B^N_{1/z}.
\end{align}
Since $\nu_*(z)\le \nu(z)$ for all $z\in Z$, it suffices to check that
\begin{align}
\label{m_nu_st1}
\nu_*(w)B_{1/w}^N \supset \cap _{z\in Z} \nu(z) B^N_{1/z}, \quad w\in Z.
\end{align}
Let $w_n\ge w$, $\nu(w_n) \underset{n\to\infty}{\to} \nu_*(w)$. Then
$$
\nu_*(w)B_{1/w}^N \supset \cap _{n\in \N} \nu(w_n) B^N _{1/w} \supset \cap _{n\in \N} \nu(w_n) B^N _{1/w_n},
$$
which implies (\ref{m_nu_st1}).

Now we set $$\nu_{**}(z) =\inf\{\nu_*(w)\cdot N^{z‐w}:\; w\le z\}.$$
Notice that
\begin{align}
\label{phi_stst_le} \nu_{**}(z)\le \nu_*(z), \quad z\in Z.
\end{align}

We show that $\nu_{**}$ is a non‐decreasing function. Indeed, let $z_1\le z_2$. If $w\le z_1$, then $\nu_*(w)\cdot N^{z_1‐w} \le \nu_*(w)\cdot N^{z_2‐w}$; hence, $\nu_{**}(z_1)\le \inf\{\nu_*(w)\cdot N^{z_2‐w}:\; w\le z_1\}$. If  $z_1\le w\le z_2$, by the monotonicity of the function $\nu_*$ we get $$\nu_*(w)\cdot N^{z_2‐w} \ge \nu_*(z_1) \ge \nu_{**}(z_1).$$
Therefore, $\nu_{**}(z_2)\ge \nu_{**}(z_1)$.

Now we prove that 
\begin{align}
\label{nn_est}
\frac{\nu_{**}(z_2)}{\nu_{**}(z_1)}\le N^{z_2‐z_1}, \quad z_1, \; z_2\in Z, \quad z_1\le z_2.
\end{align}
This yields (\ref{mon_lip}). For every $w\le z_1$ we have
$$
\nu_{**}(z_2)\le \nu_*(w)\cdot N^{z_2‐w} = \nu_*(w)\cdot N^{z_1‐w}\cdot N^{z_2‐z_1};
$$
taking infimum over $w\le z_1$, we get (\ref{nn_est}).

Now we show that
\begin{align}
\label{m_nu_st2} \cap _{z\in Z} \nu_{**}(z) B^N_{1/z} = \cap _{z\in Z} \nu_*(z) B^N_{1/z}.
\end{align}
By (\ref{phi_stst_le}), it suffices to check that for every $w\in Z$ 
\begin{align}
\label{m_nu_st3} 
\nu_{**}(w) B^N_{1/w} \supset \cap _{z\in Z} \nu_*(z) B^N_{1/z}.
\end{align}
Let $w_n\le w$, $\nu_*(w_n)\cdot N^{w‐w_n} \underset{n\to \infty}{\to} \nu_{**}(w)$. Then
$$
\nu_{**}(w) B^N_{1/w} \supset \cap _{n\in \N} \nu_*(w_n)N^{w‐w_n}B_{1/w}^N\supset \cap _{n\in \N} \nu_*(w_n)B_{1/w_n}^N.
$$
Hence, we obtain (\ref{m_nu_st3}).

Now, (\ref{m_nu_st}) and (\ref{m_nu_st2}) imply (\ref{m_cap_zz}).


\begin{Biblio}
\bibitem{vas_1} A.A. Vasil'eva, ``Kolmogorov widths of weighted Sobolev classes on a multi-dimensional domain with conditions on the derivatives of
order $r$ and zero'', arXiv:2004.06013.


\bibitem{vas_finite} A.A. Vasil'eva, ``Kolmogorov widths of the intersection of two
finite-dimensional balls'', arXiv:2007.04894v1.

\bibitem{pietsch1} A. Pietsch, ``$s$-numbers of operators in Banach space'', {\it Studia Math.},
{\bf 51} (1974), 201--223.

\bibitem{stesin} M.I. Stesin, ``Aleksandrov diameters of finite-dimensional sets
and of classes of smooth functions'', {\it Dokl. Akad. Nauk SSSR},
{\bf 220}:6 (1975), 1278--1281 [Soviet Math. Dokl.].

\bibitem{k_p_s} A.N. Kolmogorov, A. A. Petrov, Yu. M. Smirnov, ``A formula of Gauss in the theory of the method of least squares'', {\it Izvestiya Akad. Nauk SSSR. Ser. Mat.} {\bf 11} (1947), 561‐‐566 (in Russian). 

\bibitem{stech_poper} S. B. Stechkin, ``On the best approximations of given classes of functions by arbitrary polynomials'', {\it Uspekhi Mat. Nauk, 9:1(59)} (1954) 133‐‐134 (in Russian).

\bibitem{kashin_oct} B.S. Kashin, ``The diameters of octahedra'', {\it Usp. Mat. Nauk} {\bf 30}:4 (1975), 251‐‐252 (in Russian).

\bibitem{bib_kashin} B.S. Kashin, ``The widths of certain finite-dimensional
sets and classes of smooth functions'', {\it Math. USSR-Izv.},
{\bf 11}:2 (1977), 317--333.


\bibitem{gluskin1} E.D. Gluskin, ``On some finite-dimensional problems of the theory of diameters'', {\it Vestn. Leningr. Univ.}, {\bf 13}:3 (1981), 5--10 (in Russian).

\bibitem{bib_gluskin} E.D. Gluskin, ``Norms of random matrices and diameters
of finite-dimensional sets'', {\it Math. USSR-Sb.}, {\bf 48}:1
(1984), 173--182.

\bibitem{garn_glus} A.Yu. Garnaev and E.D. Gluskin, ``On widths of the Euclidean ball'', {\it Dokl.Akad. Nauk SSSR}, {bf 277}:5 (1984), 1048--1052 [Sov. Math. Dokl. 30 (1984), 200--204]

\bibitem{itogi_nt} V.M. Tikhomirov, ``Theory of approximations''. In: {\it Current problems in
mathematics. Fundamental directions.} vol. 14. ({\it Itogi Nauki i
Tekhniki}) (Akad. Nauk SSSR, Vsesoyuz. Inst. Nauchn. i Tekhn.
Inform., Moscow, 1987), pp. 103--260 [Encycl. Math. Sci. vol. 14,
1990, pp. 93--243].

\bibitem{galeev1} E.M.~Galeev, ``The Kolmogorov diameter of the intersection of classes of periodic
functions and of finite-dimensional sets'', {\it Math. Notes},
{\bf 29}:5 (1981), 382--388.

\bibitem{galeev4} E.M. Galeev,  ``Widths of functional classes and finite‐dimensional sets'', {\it Vladikavkaz. Mat. Zh.}, {\bf 13}:2 (2011), 3‐‐14.

\bibitem{bib_glus_3} E. D. Gluskin, ``Intersections of a cube with an octahedron are poorly approximated by subspaces of small dimension'', Approximation of functions by special classes of operators, Interuniv. Collect. Sci. Works, Vologda, 1987, pp. 35‐‐41 (in Russian).

\bibitem{mal_rjut} Yu. V. Malykhin, K. S. Ryutin, ``The Product of Octahedra is Badly Approximated in the $l_{2,1}$‐Metric'', {\it Math. Notes}, {\bf 101}:1 (2017), 94‐‐99.

\bibitem{kon_m_r} S. V. Konyagin, Yu. V. Malykhin, C. S. Rjutin, ``On Exact Recovery of Sparse Vectors from Linear Measurements'', {\it Math. Notes}, {\bf 94}:1 (2013),  107‐‐114.

\end{Biblio}
\end{document}